\documentclass[11pt,a4paper]{article}
\usepackage[left=25mm,right=25mm,top=24mm,bottom=23mm]{geometry}
\usepackage[T1]{fontenc}
\usepackage{lmodern,microtype,amsmath,amssymb,enumitem,tikz,fancyhdr}
\usepackage[bookmarksnumbered,colorlinks=true,linkcolor=teal!48!black,citecolor=teal!48!black,urlcolor=teal!48!black]{hyperref}
\hypersetup{pdftitle={Sparse Additive Spanners and the Cost Spectra of Group Actions},pdfauthor={Itai Benjamini},pdfsubject={Exposition on Cayley graphs and measured group actions}}
\definecolor{ink}{RGB}{25,55,61}\definecolor{mesh}{RGB}{22,93,97}
\renewcommand{\sectionmark}[1]{}
\newcommand{\Z}{\mathbb Z}\newcommand{\E}{\mathbb E}\newcommand{\Cost}{\operatorname{Cost}}
\newcommand{\Cadd}{C_{\mathrm{add}}}\newcommand{\Spec}{\mathcal S}
\newcommand{\lead}[1]{\par\smallskip\noindent\textbf{#1}\enspace}
\newcommand{\result}{\lead}
\newcommand{\parthead}[2]{\phantomsection\addcontentsline{toc}{part}{Part #1. #2}\par\medskip
\noindent{\Large\bfseries\color{ink} Part #1.\enspace #2}\par\smallskip
\markboth{#1.\enspace #2}{}}
\begin{document}
\thispagestyle{plain}
\begin{center}
{\LARGE\bfseries\color{ink} Sparse Additive Spanners and the\\[3pt]Cost Spectra of Group Actions\par}
\vspace{7pt}{\large Itai Benjamini\par}
\end{center}
\begin{abstract}\noindent
How sparse can a Cayley graph become if we retain all its vertices and allow distances to increase by at most a fixed additive constant? We study invariant random spanning subgraphs, measuring sparsity by the expected degree at a vertex. On the square grid, a random translate of a sparse mesh brings the mean degree arbitrarily close to two. If the choice must be made equivariantly from independent vertex labels, every edge is forced. In this expository note we follow these two answers through other groups. Geometry alone forces a positive gap on each non-elementary hyperbolic Cayley graph, while on one four-regular nilpotent Cayley graph the choice of action realizes every additive cost in $[1,2]$. The minimum of the action spectrum recovers the geometric cost.
\end{abstract}
\parthead{I}{Sparse additive spanners}
\section{How few edges can serve the whole grid?}\label{sec:grid}
Let $\Gamma$ be an infinite group with identity $e$ and finite generating set $S$, closed under inverses and excluding $e$. In its Cayley graph $G$ a step from $g$ to $gs$, $s\in S$, has length one. Write $d_G=d_S$ for shortest-path distance. This is the word metric for $S$. A spanning subgraph $H$ keeps every vertex. It is an \emph{additive $L$-spanner} if
\[
d_G(g,h)\le d_H(g,h)\le d_G(g,h)+L\qquad(g,h\in\Gamma).
\]
Disconnected pairs have infinite distance, so this condition forces connectivity. For a random spanner, one deterministic finite $L$ must work for all pairs almost surely.

On the square grid keep every vertical edge and the horizontal edges on every $m$th row. A uniform row offset modulo $m$ gives a translation-invariant law. Do the horizontal part of a trip on the retained row nearest its starting point. This adds at most $2\lfloor m/2\rfloor$ steps. The origin always has two vertical neighbors and has two horizontal ones with probability $1/m$. Its mean degree is $2+2/m$.

We can do better at the same order of detour. Keep every $m$th row and every $m$th column. Each vertex outside this mesh contributes one vertical edge toward its nearest mesh row, breaking ties downward. A period cell has $m^2$ vertices and $2m+(m-1)^2=m^2+1$ edges. A uniform translate therefore has mean degree $2+2/m^2$. Every vertex reaches a mesh intersection within $m$ steps. Two intersections are joined at their exact grid distance. The extra length of any trip is at most $4m$.

\begin{center}
\begin{tikzpicture}[x=4mm,y=4mm]
\begin{scope}
\foreach \x in {0,...,12} \draw[mesh,line width=.78pt] (\x,0)--(\x,6);
\foreach \y in {0,3,6} \draw[mesh,line width=1.05pt] (0,\y)--(12,\y);
\foreach \x in {0,...,12}\foreach \y in {0,...,6} \fill[ink] (\x,\y) circle[radius=1.12pt];
\end{scope}
\begin{scope}[xshift=66mm]
\foreach \x in {0,3,6,9,12} \draw[mesh,line width=1.05pt] (\x,0)--(\x,6);
\foreach \y in {0,3,6} \draw[mesh,line width=1.05pt] (0,\y)--(12,\y);
\foreach \x in {1,2,4,5,7,8,10,11} \draw[black!78,line width=.72pt]
(\x,0)--(\x,1) (\x,2)--(\x,4) (\x,5)--(\x,6);
\foreach \x in {0,...,12}\foreach \y in {0,...,6} \fill[ink] (\x,\y) circle[radius=1.12pt];
\end{scope}
\end{tikzpicture}
\par\smallskip{\small\textbf{Figure 1.} Left, every column and every third row. Right, a mesh with one attachment per outside vertex.}\end{center}

On a finite graph one counts edges directly. On an infinite Cayley graph every vertex has the same expected degree under an invariant law. This lets us ask how close the edge density can come to one while all distances are preserved to bounded error.

\newpage
\section{Counting edges in an infinite graph}\label{sec:cost}
An invariant law assigns the same distribution to $H$ and to each of its translates. The root degree $\deg_H(e)$ therefore has the same distribution at every vertex. We measure density by half its mean, counting the two ends of an undirected edge. A periodic subgraph is preserved by a finite-index subgroup of translations. Choosing a uniform translate from its finite orbit makes the law invariant.

For a fixed $L$, let $s_G(L)$ be the infimum of $\E\deg_H(e)-2$ over invariant additive $L$-spanners. Every invariant connected spanning graph on an infinite Cayley graph has mean degree at least two \cite{BLPS}. We write
\[
\Cadd(G)=\inf_{L<\infty}\ \inf_{H\ \mathrm{invariant}\ L\text{-spanner}}\tfrac12\E\deg_H(e)
=1+\tfrac12\inf_{L<\infty}s_G(L).
\]
Thus $\Cadd(G)=1$ means that for each $\epsilon>0$ some rule with a finite error has mean degree at most $2+\epsilon$. The allowed error can grow as $\epsilon$ decreases. One may equivalently require the spanner law to be invariant under all graph automorphisms. Indeed, $\Gamma$ acts simply transitively and is a uniform lattice in $\operatorname{Aut}(G)$. Averaging a translation-invariant law over $\operatorname{Aut}(G)/\Gamma$ preserves the error and the mean root degree \cite{BLPS,main}.

\lead{Why the mesh has the right order.}
For a square of side $m=\lfloor L\rfloor+1$, the sum of the two diagonal distances in $H$ is at least $4m$. Either sum of opposite side lengths is at most $2m+2L<4m$. In a tree, the two largest of these three sums agree. Thus the six shortest retained paths between the corners contain a cycle in a ball of radius $5m$. Place disjoint such balls at the points of a uniform translate of $12m\Z^2$. In each ball delete one edge of one retained cycle, using independent continuous labels to choose. The other edges of each chosen cycle survive, so connectivity survives. The resulting mean degree is still at least two. The original mean had to pay for one deletion per $(12m)^2$ vertices. Together with the mesh construction this gives $s_{\Z^2}(L)\asymp L^{-2}$ for large $L$ \cite{main}.

\lead{The cost of connectivity.}
Levitt introduced the measured cost of an equivalence relation, and Gaboriau developed its theory \cite{Levitt,Gaboriau}. For the orbit relation of a free probability-preserving action, it is the infimum of half the mean degree over measurable graphs whose connected components are those orbits. Edges may connect any two points of an orbit. The ordinary \emph{group cost} $\Cost(\Gamma)$ takes the infimum over such actions. Equivalently it is the infimum over translation-invariant random connected graphs on vertex set $\Gamma$, with unrestricted edges \cite{HP}. In particular
\[
\Cost(\Gamma)\le C_*(\Gamma):=\inf_S \Cadd(\operatorname{Cay}(\Gamma,S))
\le\Cadd(\operatorname{Cay}(\Gamma,S)).
\]
The middle quantity allows a choice of generators. In Part II we fix the generators and prescribe the action that must choose the edges.

\section{Products and polynomial growth}
\result{Cartesian products.}
For infinite Cayley graphs $X,Y$ we have $\Cadd(X\square Y)=1$. In each factor join vertices at distance between $1$ and $2r$ in an auxiliary graph of maximum degree $D$. Give each vertex an independent uniform priority in $[0,1]$. Accept it if none of its lower-priority neighbors is accepted. From a fixed vertex the chance of a decreasing path of length $n$ is at most $D^n/(n+1)!$. Thus the recursion almost surely terminates and commutes with translations. Accepted vertices are more than $2r$ apart, and every rejected vertex lies within $2r$ of an accepted one.

Choose $A\subset X$ and $B\subset Y$ independently, with root inclusion probabilities $p,q$. Each radius-$r$ ball contains at most one selected point, so $p\le |B_X(e,r)|^{-1}$ and similarly for $q$. Keep the rows $X\times\{b\}$ for $b\in B$ and columns $\{a\}\times Y$ for $a\in A$. Call their union $K$. Each vertex outside $K$ sends one edge decreasing its distance to $K$, breaking ties by generator order. A cycle of $k$ attachment edges would need $k$ distinct tails on the cycle, but a vertex nearest $K$ cannot send an edge along it.

Each attachment path reaches $K$ in at most $2r$ steps, since the selected points form $2r$-nets. Once on a retained row or column, at most $2r$ further steps along it reach $A\times B$. Between intersections $(a,b)$ and $(a',b')$, follow row $b$ to $(a',b)$ and column $a'$ to $(a',b')$. The length is $d_X(a,a')+d_Y(b,b')$, exactly the product distance. Each endpoint is within $4r$ of its chosen intersection in the spanner, so the triangle inequality gives additive error at most $16r$. For factor degrees $d_X,d_Y$,
\[
\E\deg_H(e)=d_Xq+d_Yp+2(1-p)(1-q)\ \longrightarrow\ 2.
\]
The last term counts both ends of the edge sent by each vertex outside $K$.

\lead{Nilpotent groups.}
A group is virtually nilpotent of class at most two if it contains a finite-index subgroup in which all commutators are central. Every Cayley graph of an infinite group of this kind admits \emph{periodic} additive spanners with mean degree tending to two. For a virtually $\Z^d$ group with $d\ge2$ the optimal excess has order $L^{-d}$, and for the standard integer Heisenberg graph it has order $L^{-4}$ \cite[Section~5]{main}.

The issue in a noncommutative group is how to return from a long reserved path. In the integer Heisenberg group let $z=[a,b]$, which is central. To move from $g$ to $ga^n$, transfer to a reserved $a$-line at $gz^j$, follow it to $ga^nz^j$, and undo $z^j$. A periodic construction bounds $j$ uniformly. Repeating such a transfer at every letter of a word would lose the additive bound. Bodart's lemma supplies a representative of near-minimal length consisting of a bounded number of powers, so only boundedly many transfers are needed \cite{Bodart,main}.

For a group of polynomial growth with ball volumes comparable to $r^D$, $D>1$, a weaker distance requirement is available on every free action. For $(D+1)/(2D)<\theta<1$ one can achieve mean degree arbitrarily close to two with distance at most $R+O((1+R)^\theta)$ for ambient distance $R$. Sparse nets at successive scales cover the long trips \cite{main}. Whether one can always replace this sublinear error by a constant is the first open problem below.

Class three gives a concrete obstruction. Let $F=F(a,b)$ and quotient by $\gamma_4F$, the fourth term of its lower central series. The resulting \emph{Cartan group} is the free class-three nilpotent group on two generators. On its standard graph every periodic finite-error spanner keeps every edge. The same conclusion holds for a translation-invariant law ergodic under the single translation $ab$ \cite[Thm.~5.5]{main}. Arbitrary invariant laws are still open. We will meet a different class-three group in Section~\ref{sec:engel}.

\section{Hyperbolic branching}\label{sec:hyper}
For each non-elementary hyperbolic Cayley graph $G$ there is $\eta(G)>0$ such that \emph{every} invariant finite-error spanner has mean degree at least $2+\eta(G)$ \cite[Thm.~2.1]{main}. Hyperbolicity means that geodesic triangles are uniformly thin. Non-elementary excludes finite and virtually cyclic groups. The lower bound is independent of the size of the permitted finite error.

\lead{Why lines survive.}
Take shortest $H$-paths from the identity to successive vertices of an ambient geodesic ray. A subsequence converges on every finite initial segment. At time $t$ the deficit $t-d_G(e,v_t)$ is a nondecreasing integer at most $L$, hence eventually constant. The tail is an ambient geodesic, and hyperbolicity puts it at the same boundary point as the original ray. For the attracting endpoint of an infinite-order element $g$, the retained tail stays a bounded distance from $\{g^n:n\ge0\}$. By recurrence there are translations $g^{-n_j}$ that return the edge configuration on larger and larger balls about $e$. Each translated tail passes through a common finite ball. A subsequential limit in the original $H$ is a full geodesic line joining the two ends of $g$. These pairs of fixed points are dense among pairs of distinct boundary points, and compactness supplies a retained line between each pair of distinct boundary points.

\lead{Why the mean degree must rise.}
Take the union $K$ of all retained ambient lines. Fix disjoint boundary neighborhoods so that every line between them crosses one finite ball $B$. Choose more than $|B|$ pairs with all endpoints different in those neighborhoods. Two of the corresponding lines must meet. Compactness shows that they cannot all avoid branching in arbitrarily large balls: a limit of maximum degree two would make two intersecting lines coincide despite their distinct endpoints. Consequently every translate of some fixed ball $B_R$ contains a branch of $K$.

Each vertex outside $K$ chooses one retained edge decreasing its distance to $K$. These edges are distinct, and their expected incoming and outgoing counts agree by translation invariance. With $\deg_K(e)=0$ outside $K$ this gives the useful attachment bound
\[
\E\deg_H(e)\ \ge\ \E\deg_K(e)+2\mathbb P(e\notin V(K)).
\]
Vertices of $K$ have degree at least two and a proportion at least $1/|B_R|$ are branches. Thus $\eta(G)$ can be taken to be $1/|B_R|$. This proof does not control $R$ as the generators vary on a fixed group.

\lead{Rigidity of a tiling.}
The argument has a sharper counterpart for transitive tilings. For the regular hyperbolic tiling $\{p,q\}$, with $p$-gonal faces and degree $q$, $1/p+1/q<1/2$, every finite-error spanner invariant under its full symmetry group keeps every edge almost surely \cite[Thm.~3.1]{main}. In particular the $7$-regular triangulation admits no thinning of this kind. Choose two incident edges whose midpoints are maximally far apart. Successive half-turns extend them to a periodic path whose midpoints lie on one hyperbolic line. Projection along that line assigns the same maximal progress to each forward step of the path. All other midpoint transitions lose a uniform positive amount. A route with bounded extra length can therefore bypass only boundedly many of the prescribed edges, however long a segment it follows. Stationarity along the period forces zero probability of omission. Tiling symmetries carry this conclusion to every edge.

\begin{center}
\begin{tikzpicture}[x=1.5cm,y=3.75cm,font=\small]
\draw[mesh,dashed] (-2.62,0)--(2.68,0) node[right,black] {$\lambda$};
\draw[ink,line width=1pt] plot coordinates {(-2.50000,-0.11104) (-2.41667,-0.09110) (-2.33333,-0.07196) (-2.25000,-0.05344) (-2.16667,-0.03537) (-2.08333,-0.01761) (-2.00000,0.00000) (-1.91667,0.01761) (-1.83333,0.03537) (-1.75000,0.05344) (-1.66667,0.07196) (-1.58333,0.09110) (-1.50000,0.11104) (-1.50000,0.11104) (-1.41667,0.09110) (-1.33333,0.07196) (-1.25000,0.05344) (-1.16667,0.03537) (-1.08333,0.01761) (-1.00000,-0.00000) (-0.91667,-0.01761) (-0.83333,-0.03537) (-0.75000,-0.05344) (-0.66667,-0.07196) (-0.58333,-0.09110) (-0.50000,-0.11104) (-0.50000,-0.11104) (-0.41667,-0.09110) (-0.33333,-0.07196) (-0.25000,-0.05344) (-0.16667,-0.03537) (-0.08333,-0.01761) (0.00000,0.00000) (0.08333,0.01761) (0.16667,0.03537) (0.25000,0.05344) (0.33333,0.07196) (0.41667,0.09110) (0.50000,0.11104) (0.50000,0.11104) (0.58333,0.09110) (0.66667,0.07196) (0.75000,0.05344) (0.83333,0.03537) (0.91667,0.01761) (1.00000,-0.00000) (1.08333,-0.01761) (1.16667,-0.03537) (1.25000,-0.05344) (1.33333,-0.07196) (1.41667,-0.09110) (1.50000,-0.11104) (1.50000,-0.11104) (1.58333,-0.09110) (1.66667,-0.07196) (1.75000,-0.05344) (1.83333,-0.03537) (1.91667,-0.01761) (2.00000,0.00000) (2.08333,0.01761) (2.16667,0.03537) (2.25000,0.05344) (2.33333,0.07196) (2.41667,0.09110) (2.50000,0.11104)};
\fill (-2.500,-0.11104) circle[radius=1.35pt];
\fill (-1.500,0.11104) circle[radius=1.35pt];
\fill (-0.500,-0.11104) circle[radius=1.35pt];
\fill (0.500,0.11104) circle[radius=1.35pt];
\fill (1.500,-0.11104) circle[radius=1.35pt];
\fill (2.500,0.11104) circle[radius=1.35pt];
\draw[black,fill=white,line width=.65pt] (-2,0) circle[radius=2pt];
\draw[black,fill=white,line width=.65pt] (-1,0) circle[radius=2pt];
\draw[black,fill=white,line width=.65pt] (0,0) circle[radius=2pt];
\draw[black,fill=white,line width=.65pt] (1,0) circle[radius=2pt];
\draw[black,fill=white,line width=.65pt] (2,0) circle[radius=2pt];
\end{tikzpicture}
\par\smallskip{\small\textbf{Figure 2.} In the seven-regular triangulation the vertices turn, while successive edge midpoints (open circles) advance equally far on the hyperbolic line $\lambda$. The coordinates are distance along the line and signed distance from it.}
\end{center}

\section{What survives a change of generators?}\label{sec:change}
The graph $\operatorname{Cay}(\Gamma,S)$ specifies which moves have length one. Even adding a few types of move can force a spanner to keep more edges.

\lead{Two products of the same trees.}
On $F_2\times F_2$ the Cartesian generators give cost one. Now also allow simultaneous generator steps in both coordinates. The ambient distance becomes the maximum of the two tree distances. On a pair of rays, a nearly shortest route eventually has to advance in both coordinates at every step. The full argument forces all sixteen diagonal edges on one parity class of $|x|+|y|\bmod 2$. Their contribution to the mean degree is eight, and the attachment bound adds one for the other class. This is sharp: retain the diagonals on one uniformly chosen parity class and attach each remaining vertex once. Two-step backtracks in the shorter coordinate realize the maximum distance within the chosen class. The mean degree is nine and the additive error is at most four \cite[Thm.~6.1]{main}.

\lead{A tree of triangles.}
With the standard generators $a^{\pm1},b^{\pm1}$ of $C_3*C_3$, independently deleting a uniform side in every triangle leaves a connected graph of cost $4/3$, the ordinary group cost. A route from $e$ to $(ab)^n$ crosses $2n$ triangles at compulsory entry and exit vertices. Each omitted direct side adds a step. Translation invariance makes the expected loss grow linearly if either generator has positive omission probability. Every additive spanner therefore has cost two.

A different generating set gives $5/3$. With $A=\langle a\rangle$, use $(A\setminus\{e\})\cup AbA\cup Ab^{-1}A$. The adjacent three-point fibers of the six-regular quotient graph are completely joined. Choose one representative per fiber periodically, copy the quotient edges on the representatives and attach the other two vertices. The additive error is at most two and the cost is $(6+4)/6=5/3$. For \emph{each} fixed finite generating set the cost is strictly above $4/3$ \cite[Thms.~2.11--2.12]{main}. The proof passes a hypothetical sequence approaching ordinary cost to a forest carrying all boundary lines. Finite local rules then give a periodic additive spanning tree. Rigidity of tree actions would put all its branch points in a single orbit, inconsistent with the stabilizers of the two distinct order-three factors. Whether the infimum over generating sets is $4/3$ remains open.

\lead{An infimum approached by changing generators.}
For the binary lamplighter $\bigoplus_{\Z} C_2\rtimes\Z$, let $t$ move the cursor one site, and $a$ toggle the lamp there. Every fixed finite generating set has additive cost strictly above one. In a hypothetical limit of mean degree two, the components are trees with at most two ends. Their double rays give only countably many possible limiting lamp configurations. Bounded detours along independently prescribed switches would force a match with a nonatomic random configuration \cite[Thm.~7.1]{main}. Yet for $a_r=t^rat^{-r}$, $r\ge1$, the four moves $t,at,a_rt,aa_rt$ and their inverses have least mean degree $2+2^{1-r}$, attained with error at most $4r$ \cite[Thm.~7.3]{main}. Thus $C_*=1$ for this group.

\lead{A gap for every generating set at once.}
For a closed orientable hyperbolic three-manifold $M$, write $\Gamma=\pi_1M$ and $b=\dim H_1(M;\mathbb R)$. Then \cite[Thm.~1.3]{main}
\[
\Cost(\Gamma)=1,\qquad C_*(\Gamma)\ge 1+\lfloor b/2\rfloor\quad (b\ge4).
\]
For connectivity alone, a finite cover fibers over the circle. Keep the height-changing edges and place all fiber generators at sparse periodic levels. This proves cost one in the cover, and the ordinary cost calculation passes to the finite extension.

To see why the metric costs more, project $S$ to real homology and take its convex hull. Exposing one of its edges by a linear functional picks out generators that advance at maximal speed. A long route in those directions has a bounded budget for slower steps. For $b\ge4$ their subgroup spans a plane in homology, too large to be virtually cyclic and too small to be a virtual fiber or a finite-index subgroup. Three-manifold geometry gives a quasiconvex subgroup whose commensurator has ordinary cost at least two. Here the commensurator consists of elements conjugating the subgroup to one with finite-index intersection in each. The proof counts at least one unit above the connectivity baseline for each of $\lfloor b/2\rfloor$ disjoint generator families selected on opposite facets of the polytope. One can find this many families on a facet of the $b$-dimensional homology polytope. Its graph has minimum degree at least $b-1$, hence a path through at least $b$ vertices. Take alternate edges on that path. Their negatives belong to the opposite facet. The selected families use disjoint edge labels, so their contributions add \cite[Prop.~2.8]{main}. Finite covers can have arbitrarily large $b$ \cite{Agol}. This gives an unbounded gap above ordinary cost even after optimizing the generators.

A further comparison explains why the open hyperbolic question fixes the group. For $F_2\times C_m$, $m\ge3$, the optimized cost is exactly $1+1/m$ \cite[Section~2]{main}. A spanning path in each cycle fiber and one uniformly chosen tree layer give the upper bound. The index-$m$ free subgroup gives the matching lower bound through $\beta_1^{(2)}=1/m$ and Gaboriau's cost inequality \cite{GabBetti}. Thus a uniform gap over \emph{all} non-elementary hyperbolic groups is impossible.

The lamplighter forces a positive excess for each generating set separately, yet the excess tends to zero along a sequence of generating sets. For the three-manifold groups with $b\ge4$, one lower bound works for all generating sets. We turn next to a different variation: what can change without moving a single edge of the ambient Cayley graph?

\clearpage
\parthead{II}{The spectra of group actions}
\section{Who chooses the edges?}\label{sec:action}
The first construction on the grid used a common row offset. Its law is translation invariant, but a given probability space with its translations may have no way to supply such an offset. We now keep the group, the generators and the word distance fixed. We ask what spanners can be chosen from a prescribed action, using no randomness beyond the state of that action.

\lead{The action supplies the randomness.}
Draw a state $x$ from a standard probability space $(X,\mu)$. A right action gives a measurable map $x\mapsto xg$ for each $g$, with $xe=x$ and $(xg)h=x(gh)$. Think of $xg$ as the same state viewed from vertex $g$. We assume every map $x\mapsto xg$ preserves $\mu$. The action is assumed to be \emph{free}, meaning that almost surely $xg=x$ only for $g=e$. Thus the orbit $\{xg:g\in\Gamma\}$ is a copy of $G$ when we join $xg$ to $xgs$. These graphs on the orbits together form the \emph{Cayley graphing}.

\lead{An edge rule.}
For each $s\in S$ choose a measurable event $A_s\subseteq X$. In state $x$, form a spanning subgraph $H_x$ of $G$ by keeping $\{g,gs\}$ precisely when $xg\in A_s$. Require $y\in A_s\iff ys\in A_{s^{-1}}$ almost surely, so both endpoints agree. This describes every measurable choice of edges in the Cayley graphing. Changing the viewpoint from $e$ to $g$ translates the whole output graph by $g^{-1}$. This is what it means for the rule to be equivariant.

\lead{Admissibility and cost.}
For a nonnegative integer $L$, let $\mathcal A_L(X,S)$ be the families $A=(A_s)_{s\in S}$ above for which, for $\mu$-almost every $x$ and \emph{every} $g,h\in\Gamma$,
\[
d_{H_x}(g,h)\le d_S(g,h)+L.
\]
We call $H_x$ an \emph{additive $L$-spanner}. One $L$ must work for all pairs and almost all states. At the identity the degree is $\deg_{H_x}(e)=\sum_{s\in S}\mathbf1_{A_s}(x)$. The degree counts an undirected edge at both ends. We therefore define the \emph{cost of the rule} as $\tfrac12\int_X\deg_{H_x}(e)\,d\mu(x)$. The \emph{additive cost of the action} is
\[
\boxed{\displaystyle
\Cadd(X,S)=\inf_{L=0,1,\ldots}\ \inf_{A\in\mathcal A_L(X,S)}
\frac12\sum_{s\in S}\mu(A_s).}
\]
Here $X$ includes its action. The full graph qualifies with $L=0$ and cost $|S|/2$. Ordinary measured cost on these orbits is a lower bound since every admissible rule connects them.

An action is \emph{ergodic} if every event invariant under all the group elements has probability zero or one. Keeping $\Gamma$ and $S$ fixed, we define the \emph{ergodic additive-cost spectrum} by
\[
\Spec(\Gamma,S)=
\big\{\Cadd(X,S):\ \Gamma\curvearrowright(X,\mu)
\text{ is free, ergodic and probability-preserving}\big\}.
\]
A spectrum records which costs are attained by actions. For any one action its cost is an infimum and may be approached by a sequence of spanners.

\lead{The link between the two costs.}
For every finite symmetric generating set of an infinite group,
\[
\boxed{\displaystyle\Cadd(\operatorname{Cay}(\Gamma,S))
=\min\Spec(\Gamma,S).}
\]
Every action rule gives an invariant law. For the converse choose a translation-ergodic spanner law $\nu_n$ within $2^{-n}$ of the geometric cost for each $n$. Such a law exists by decomposing a near-minimizer into ergodic components: the error bound holds on almost all of them and their costs average to the original cost.

Take the product of the $\nu_n$ and an independent Bernoulli shift (independent uniform labels on $\Gamma$). Almost every ergodic component for the diagonal action has every prescribed marginal, since each marginal is ergodic and there are only countably many. On one such component the labels ensure freeness, and reading the $n$th graph gives a rule of cost tending to $\Cadd(G)$. Thus one action attains the minimum even when none of its individual rules attains its cost.

\section{The square grid revisited}

Return to $S=\{\pm u,\pm v\}$, with $u=(1,0)$ and $v=(0,1)$. The row construction of Figure~1 uses only a single shared residue. We can now ask precisely which actions supply it.

\lead{A shared row offset.}
Fix $k\ge1$. The state consists of independent uniform $[0,1]$ labels on $\Z^2$ and an independent uniform residue $r$ modulo $k$. Translation shifts the labels. Moving north also increases $r$ by one, while moving east leaves $r$ fixed. This action is free and ergodic. Take
\[
A_v=A_{-v}=X,\qquad A_u=A_{-u}=\{r=0\}.
\]
Thus all vertical edges are kept, and horizontal edges are kept on the rows $j$ with $r+j=0\bmod k$. Its cost is therefore
\[
\tfrac12(1+1+1/k+1/k)=1+1/k .
\]

Section~\ref{sec:grid} gives the error bound $2\lfloor k/2\rfloor$. Thus $\Cadd(X,S)\le1+1/k$ for this action. To prove equality we need to rule out all cheaper rules on the same probability space.

\lead{Independent labels alone.}
Drop the residue and let $X=[0,1]^{\Z^2}$ with independent uniform coordinates. The \emph{Bernoulli action} translates the whole array of independent labels. The horizontal shift is ergodic, since sufficiently distant translates of events depending on finitely many labels are independent. Suppose a rule omits a horizontal edge with positive probability. Ergodicity implies that almost surely each row has omitted edges arbitrarily far to the east. There are only countably many rows.

Yet a path from $(0,0)$ to $(n,0)$ of length at most $n+L$ has at most $L$ steps other than east. Choose one such path for each $n$. The grid has finite degree, so a subsequence agrees on each initial segment and gives an infinite path with at most $L$ steps other than east. Its tail runs east forever on a single row. This contradicts the omitted edges there. The vertical argument is identical. The full grid is therefore forced, giving $\Cadd(X,S)=2$. The shared offset matters because it supplies an event of probability $1/k$ that is unchanged on moving east.

\lead{The exact square-grid formula.}
For the horizontal shift $u$, let $\mathcal I_u$ consist of the measurable events that are unchanged by $u$, up to null sets. An \emph{atom} of $\mathcal I_u$ is an event of positive probability containing no event of this algebra of strictly intermediate probability. In an ergodic $\Z^2$-action, if there is one atom, its translates are disjoint or equal up to null sets, have equal mass, and cover almost all of $X$. Write $m_u$ for their number, and $m_u=\infty$ when there are no atoms. Define $m_v$ the same way. Then
\[
\Cadd(X,S)=1+\frac{1}{m_um_v},
\qquad
\Spec(\Z^2,S)=\{1\}\cup\{1+1/k:k=1,2,\ldots\},
\qquad 1/\infty=0.
\]
For the action with a row offset, the horizontal shift leaves precisely the residue as invariant information, so $m_u=k$. The $k$th power of the vertical shift is ergodic within each residue and one vertical step cycles through them, so $m_v=1$. Thus the explicit row-offset rule has the least possible cost.

\lead{Where the product comes from.}
Write $R$ for the event that the row through the origin has every edge, and $D$ for the corresponding column event. The eastward-tail argument still applies. By recurrence under the horizontal shift, almost surely every omitted edge has omitted edges arbitrarily far ahead on its row. Thus an eastward tail lies on an intact row. Almost every configuration has such a row somewhere, so $\mu(R)>0$ by translation invariance and countability. Likewise $\mu(D)>0$.

The intact lines contribute $\mu(R)+\mu(D)$ to the cost. Every vertex outside them can choose one edge along a shortest retained path toward them. These edges are distinct, and their mean indegree equals their mean outdegree by invariance. The extra cost is $1-\mu(R\cup D)$, giving the lower bound $1+\mu(R\cap D)$. Now $R\in\mathcal I_u$ and $D\in\mathcal I_v$. The shift $u$ is ergodic on $\mathcal I_v$, since an event fixed by both shifts is fixed by all of $\Z^2$. Averaging the indicator of $D$ along $u$ therefore gives $\mu(D)$, while $R$ stays fixed. Hence $\mu(R\cap D)=\mu(R)\mu(D)\ge 1/(m_um_v)$.

For the reverse bound choose rows and columns from one atom in each finite algebra. They recur with bounded gaps. In an algebra with no atoms one can choose lines of arbitrarily small positive density with bounded gaps, by choosing a measurable maximal set of positions separated by more than $r$ steps of the transverse shift. Separation makes the selection sparse and maximality brings every position within $r$ steps of it. Join the other vertices to this mesh by one edge each, pointing toward the mesh. Every vertex can reach a mesh intersection in boundedly many steps, and intersections are joined at their exact grid distance. This gives the matching infimum.

Cost $1$ occurs, for example, if the action supplies compatible uniform row offsets modulo $2,4,8,\ldots$ as well as independent labels. The same action permits the preceding row construction for every $k=2^j$. Each rule has its own finite error bound. Every individual finite-error rule still costs strictly more than one, since $\mu(R)\mu(D)>0$.

\lead{Other lattice generators.}
For a finite symmetric generating set $S$ of $\Z^d$, the Bernoulli cost equals half the number of vertices of the convex hull of $S$ in $\mathbb R^d$. For a generator $s$ at a vertex of that convex hull, a linear functional strictly maximized at $s$ bounds the number of other steps in a near-shortest route to $ns$. The intact-line argument forces that direction. Conversely, rounding a convex representation routes a displacement using a bounded number of runs in these directions and a bounded remainder. Sparse connecting edges handle the remainder at arbitrarily small extra cost. The coordinate generators give cost $d$. Adding $\pm(1,1)$ to the square-grid generators gives a hexagon and raises the Bernoulli cost from two to three.

\section{An interval on a specific Engel lattice}
\label{sec:engel}

The square-grid spectrum is discrete away from its limit point $1$. For an arbitrary finite generating set of $\Z^d$ the full lattice formula still gives only rational costs. Call a group torsion-free if every nonidentity element has infinite order. For torsion-free finitely generated nilpotent groups of class at most two, the cost is one plus a minimum of overlaps of atoms in finitely many invariant algebras. The relevant finite factors have a finite common refinement with equally weighted atoms. This proves rationality for every free ergodic action \cite[Sections~3--5]{AddSpec}. Class three is the first possible class for an irrational value, and one four-regular graph already realizes the entire interval $[1,2]$.

\lead{Which Engel group?}
The real Engel group is the unique connected, simply connected Lie group with Lie algebra having basis $B,A,C,D$ and nonzero brackets $[B,A]=C$, $[B,C]=D$, up to antisymmetry \cite{EngelLie}. For the word metric we specify an integer lattice. Throughout this paper the \emph{discrete Engel group} is
\[
E=\Z^3\rtimes_U\Z,\quad Ue_0=e_0+e_1,\quad Ue_1=e_1+e_2,\quad Ue_2=e_2,\qquad
(m,i)(n,j)=(m+U^i n,i+j).
\]
Unless otherwise stated its generators are $a=(e_0,0)$, $b=(0,1)$ and their inverses. Conjugation by $b$ acts as $U$ on $\Z^3$. Starting from $e_0$, successive commutators give $e_1$ and $e_2$. Thus the group, its generators and its four-regular graph are completely specified.

The real group admits other lattices. For example $E$ is a lattice in $\mathbb R^3\rtimes_{U^t}\mathbb R$, where $U^t=\exp(t\log U)$ and $\log U$ has a single nilpotent block of size three. Its index-two subgroup $E'=\Z^3\rtimes_{U^2}\Z$ is a lattice there as well. Their abelianizations are $\Z^2$ and $\Z^2\oplus\Z/4\Z$, respectively. Indeed $(U-I)\Z^3=\langle e_1,e_2\rangle$, whereas $(U^2-I)\Z^3=\langle 2e_1+e_2,2e_2\rangle$. The differing torsion in the abelianizations distinguishes the two lattices. The free class-three Cartan group from Part I has two independent commutators of weight three, whereas $E$ has one. The following calculation concerns $E$ with the generators just specified.

\lead{Why an intact line must occur at a prescribed height.}
Write $M=\Z^3\le E$ and call the second coordinate of $(m,i)$ its height.

At height $i$, a positive $a$-step adds $U^ie_0=(1,i,\binom i2)$ to the $M$-coordinate. On a path to $a^n$ of length at most $n+L$, the first coordinate bounds the number of negative $a$-steps and $b^{\pm1}$-steps by $L$, so $|i|\le L$. The last two coordinates vanish at the endpoint. Since $i^2=i+2\binom i2$, the sum of $i^2$ over positive $a$-steps equals its sum over negative ones. The latter is bounded by $L^3$. All but a bounded number of steps therefore go forward in direction $a$ at height zero. As on the grid, a limiting path and recurrence force an intact $a$-line. Here it must lie in the starting $M$-orbit.

\lead{Where the continuum enters.}
Let $\pi_a(x)$ record all information unchanged by $a$. Formally, events determined by $\pi_a$ are exactly the $a$-invariant measurable events, modulo null sets. Because $M$ commutes with $a$, it acts on these records. Write $m(x)$ for the number of distinct records $\pi_a(xm)$ as $m$ ranges over $M$. If the transformation $x\mapsto xb$ is ergodic, the height coordinate forces intact $b$-lines by the limiting-ray argument of Section~\ref{sec:hyper}. On each finite orbit of records an intact $a$-line must be selected at least once. Its minimum relative density there is $1/m(x)$. Infinite orbits admit arbitrarily sparse selections within bounded reach. The resulting formula is
\[
\Cadd(X,S)=1+\int_X\frac{1}{m(x)}\,d\mu(x).
\]
For the upper bound, travel to a selected $a$-line by an element of $M$, follow it for $n$ steps, and return by the inverse element. Commutativity gives the correct endpoint. The proof supplies bounded detours for these transfers at arbitrarily small extra cost and writes any trip as a bounded number of $a$- and $b$-runs within bounded error of shortest length.

To realize $1+p$, choose rows $i\in\Z$ independently. With probability $p$ fill a row with stars, and otherwise fill it with independent fair bits. The element $b$ shifts the independent rows and hence acts ergodically. The element $a$ shifts row $i$ by $-i$ sites. An independent Bernoulli action of the group is included to ensure freeness. The $a$-invariant record consists of all row types and the contents of row zero. The group $M$ preserves the types and $e_1$ shifts row zero by one site. Its record thus has orbit size one on a star row and infinite orbit size almost surely on a bit row. The integral is $p$, so every $1+p$ occurs. The general lower bound is one, and the full four-regular graph gives the upper bound two on every action. Thus $[1,2]$ is the entire spectrum for these generators. The family can even be chosen so that, for each fixed $g$, its action as a single transformation is the same up to a measure-preserving change of coordinates for all $p$. The change of coordinates can depend on $g$.

These actions can be \emph{weakly mixing}, meaning that acting simultaneously on two independent copies of $X$ is ergodic. With another generating set, the same family attains $[2,3]$ even though each generator acts as a Bernoulli shift. There a word of length two detects the variation. On the four-regular graph every free action has infimum one if we allow $d_H\le d_S+f(d_S)+L$ for any prescribed nonnegative $f(n)\to\infty$.

\section{What survives strong mixing?}

An action is \emph{strongly mixing} if $\mu(A\cap Bg)\to\mu(A)\mu(B)$ for measurable $A,B$ as $g$ leaves every finite set, where $Bg=\{xg:x\in B\}$. Every infinite-order element then acts ergodically. For each generating set of a torsion-free nilpotent group of class at most two, all strongly mixing actions have the same cost. For a direct product of a free group and a finite cyclic group, one can obtain a whole interval even in this class of actions \cite[Section~9]{AddSpec}.

\lead{The same orbits, different costs.}
Let $F_{25}$ be the free group on $a_1,\ldots,a_{25}$, and let $c$ generate the cyclic group $C_3$ of order three. In $F_{25}\times C_3$ use $c^{\pm1}$ and all $a_j^{\pm1}c^{\pm1}$, with the signs chosen independently. There are free strongly mixing actions on the same probability space with exactly the same orbits whose additive costs attain every value in $[9,86/9]$. Their common orbit relation has ordinary measured cost nine.

Here each $C_3$-orbit is a triangle. Collapsing these triangles gives the $50$-regular tree. Across a tree edge, each point of one triangle can reach two of the three points in the next. In one action we can measurably select a point per triangle so that adjacent selected points are joined. Keep their tree and attach the other two points of each triangle to its selected point. This costs $(25+2)/3=9$ edges per vertex and adds at most two steps to any trip.

For independent labels, suppose two successive tree edges each have only one retained crossing. These crossings must meet. A positive probability of failure would force expected detour to grow linearly along repetitions of this two-step pattern on a tree geodesic. In each triangle choose the common endpoint of its single crossings, or any measurable representative if there are none. Color $x$ by the number of $c$-steps to this representative, modulo three. Equal colors at $x$ and $xa_j$ force a second crossing. Let $m$ be the probability of equal colors across a random free-generator step. Encode the colors by the three complex cube roots of unity. They average to zero on each triangle, and the real part of their correlation across a step is $m-(1-m)/2$. For Bernoulli labels, Kesten's norm bound \cite{Kesten} $7/25$ for averaging on the $50$-regular tree gives $m-(1-m)/2\ge-7/25$, hence $m\ge11/75$. There are $25$ unoriented tree edges per triangle on average and three vertices per triangle. The crossings alone cost at least $(25/3)(1+11/75)=86/9$.
The interval comes from varying how tree steps permute the three points above a vertex. Independently on each tree edge, choose a displacement in $C_3$, zero with probability $p$ and otherwise an equally likely choice of $+1$ or $-1$. Keep the same independent uniform and sign labels to make these choices jointly for $0\le p\le1/3$. The free group moves along the tree and adds these displacements. The element $c$ rotates the three points over each vertex. Hence the orbit of a state is the entire three-point fiber over every state in its base orbit, independently of $p$. At $p=0$ adjacent representatives can always be joined. At $p=1/3$ the displacement is uniform and the action is Bernoulli. Along a long reduced word the independent displacements mix the fiber coordinate. Changing $p$ alters only a small fraction of crossings, and repairing routes near them costs proportionally little \cite[Section~9]{AddSpec}. The cost varies continuously, so the intermediate value theorem supplies the interval.

\lead{A Kazhdan comparison.}
Let $K=\mathrm{SL}_3(\Z)$ and in $K^k$, $k\ge2$, allow the elementary matrices $I\pm E_{ij}$ in one coordinate at a time. All free actions have ordinary cost one \cite{Khezeli}. The action by translations on the profinite completion (the inverse limit of the finite quotients) has additive cost one, whereas the Bernoulli action has additive cost at least $2k$. For an elementary matrix $s$, the word $(ss^{\mathsf T})^n$ grows at the maximal possible rate in operator norm. An inefficient adjacent pair loses a definite factor, so a route with bounded extra length can have only boundedly many such pairs. A limiting ray alternates a generator and its transpose. Bernoulli ergodicity forces all four incident edges of at least one such pair in each coordinate \cite[Thm.~9.2]{main}.

\lead{How sharp is the bounded-error requirement?}
Write $C_c(X,S)$ for the action cost defined with distance bound $d_H\le c\,d_S+L$, where $c\ge1$. A finitely generated group is \emph{amenable} if it has finite nonempty sets with edge boundary divided by set size arbitrarily small. For each generating set of such a group, $\lim_{c\downarrow1} C_c(X,S)$ is the same for every free action. The proof transfers edge statistics on large finite balls between free amenable actions and joins the local routes, absorbing the joining error into the extra multiplicative allowance \cite[Section~13]{AddSpec}. In the $F_{25}\times C_3$ example even these limiting values contain the interval $[9,86/9]$.

\clearpage
\section{Open problems}\label{sec:problems}\markboth{Open problems}{}
Two choices run through these questions. Allowing every invariant law asks for the \emph{bottom} of the spectrum. Prescribing the action asks which values above it occur. Changing generators changes the distance constraint in both questions.

\subsection*{Geometry and the choice of generators}
\begin{enumerate} 
\item \textbf{Polynomial growth.} Does every infinite Cayley graph of polynomial growth have $\Cadd(G)=1$? For virtually class-two nilpotent groups of growth degree $D>1$, is the optimal fixed-error excess $s_G(L)$ of order $L^{-D}$? This rate holds for virtually $\Z^d$, $d\ge2$, and for the standard Heisenberg graph.

\item \textbf{The Cartan test.} Can an invariant finite-error spanner of the standard Cartan graph omit any edge with positive probability? Every periodic such spanner is the full graph. If full rigidity holds for all invariant laws, this graph has additive cost two and gives a negative answer to the first question.

\item \textbf{A uniform hyperbolic gap.} For a fixed non-elementary hyperbolic group $\Gamma$, is $C_*(\Gamma)>1$? Equivalently, is there an $\eta_\Gamma>0$ such that every finite-error invariant spanner of every Cayley graph of $\Gamma$ has mean degree at least $2+\eta_\Gamma$? The fixed-graph branching argument gives a constant depending on $S$. For closed orientable hyperbolic three-manifolds the uniform gap holds when $b_1\ge4$. It is open in general for $b_1\le3$, including rational homology spheres. In the positive range, when is the homological lower bound sharp?

\item \textbf{Hyperbolic groups of ordinary cost one.} Does the uniform gap hold for cocompact lattices in $\mathrm{Sp}(n,1)$, $n\ge2$? These groups have property~(T) and ordinary cost one \cite{HP}. In Gromov's density model with fixed generator number and density $1/3<d<1/2$, the random group is hyperbolic and has property~(T) with high probability \cite{Ashcroft}. Does it also satisfy $C_*>1$ with high probability?

\item \textbf{A sharp free-product value.} Is $C_*(C_3*C_3)=4/3$? Every individual generating set gives strictly more than $4/3$, and one gives $5/3$. The lamplighter shows how a positive gap at every generating set can disappear in the infimum.

\item \textbf{Intermediate growth.} Does the first Grigorchuk group, with its standard generators, admit invariant finite-error spanners with mean degree arbitrarily close to two?
\end{enumerate}
\subsection*{The spectrum for a fixed action class}
\begin{enumerate} 
\item \textbf{Amenability and strong mixing.} Can two free strongly mixing actions of one finitely generated amenable group have different additive costs for the same generating set? The question is open despite the interval for weakly mixing actions on $E$. Both actions would have the same limit as the multiplicative stretch decreases to one.

\item \textbf{One set of generators for every action.} Does every finitely generated nonabelian torsion-free nilpotent group admit a finite generating set giving additive cost one on \emph{every} free ergodic action? This holds in class two. For $E$, optimizing over generating pairs still gives a range of $[1,2]$, while one generating triple gives cost one for every free action \cite[Section~8]{AddSpec}.

\item \textbf{Torsion and irrationality.} Can a finitely generated nilpotent group of class two with torsion have an irrational additive cost for a free ergodic action? All its free strongly mixing actions have the same cost at a fixed generating set \cite[Section~10]{AddSpec}. Rationality for general ergodic actions is open.
\end{enumerate}

\lead{Further reading.}
The geometric theorems and their proofs are in \cite{main}. The action formulas and interval constructions are in \cite{AddSpec}. Both manuscripts are available from Itai Benjamini at \href{mailto:itai.benjamini@gmail.com}{itai.benjamini@gmail.com}.

\lead{Acknowledgment.} The paper was written by Chat gpt.

\phantomsection

\addcontentsline{toc}{section}{References}\begin{thebibliography}{99}
\small\setlength{\itemsep}{1.5pt}

\bibitem{main} I. Benjamini, \emph{Sparse Additive Spanners of Cayley Graphs: Invariant Random Spanners and Group Cost}, \mbox{AI note}, 24 November 2026.
\bibitem{AddSpec} I. Benjamini, \emph{Ergodic additive-cost spectra of Cayley graphings}, AI note, 2026.
\bibitem{Agol} I. Agol, \emph{Virtual properties of 3-manifolds}, Proc. Int. Congress of Mathematicians, Seoul 2014, vol.~I, 141--170, Corollary 6.8.
\bibitem{EngelLie} C. Bellettini and E. Le Donne, \emph{Regularity of sets with constant horizontal normal in the Engel group}, Comm. Anal. Geom. \textbf{21} (2013), 469--507, \S1.2.
\bibitem{BLPS} I. Benjamini, R. Lyons, Y. Peres and O. Schramm, \emph{Group-invariant percolation on graphs}, Geom. Funct. Anal. \textbf{9} (1999), 29--66.
\bibitem{Bodart} C. Bodart, \emph{Intermediate geodesic growth in virtually nilpotent groups}, Groups Geom. Dyn., published online 2025. \href{https://arxiv.org/abs/2306.10381}{arXiv:2306.10381}, Lemmas 1.7 and 2.1.
\bibitem{Gaboriau} D. Gaboriau, \emph{Co\^ut des relations d'\'equivalence et des groupes}, Invent. Math. \textbf{139} (2000), 41--98.
\bibitem{GabBetti} D. Gaboriau, \emph{Invariants $\ell^2$ de relations d'\'equivalence et de groupes}, Publ. Math. Inst. Hautes \'Etudes Sci. \textbf{95} (2002), 93--150.
\bibitem{HP} T. Hutchcroft and G. Pete, \emph{Kazhdan groups have cost 1}, Invent. Math. \textbf{221} (2020), 873--891.
\bibitem{Ashcroft} C. J. Ashcroft, \emph{Property (T) in density-type models of random groups}, Groups Geom. Dyn. \textbf{17} (2023), 1325--1356.
\bibitem{Kesten} H. Kesten, \emph{Symmetric random walks on groups}, Trans. Amer. Math. Soc. \textbf{92} (1959), 336--354.
\bibitem{Khezeli} A. Khezeli, \emph{Products of infinite countable groups have fixed price one}, \href{https://arxiv.org/abs/2509.08325}{arXiv:2509.08325}.
\bibitem{Levitt} G. Levitt, \emph{On the cost of generating an equivalence relation}, Ergodic Theory Dynam. Systems \textbf{15} (1995), 1173--1181.
\end{thebibliography}
\end{document}